\documentclass[12pt]{article}

\usepackage{amsfonts,amsmath,amssymb,latexsym,xcolor,mathrsfs,tikz,breqn,extarrows}
\usepackage{amsthm,authblk}

\usepackage{amssymb}
\usepackage{amsmath}
\usepackage{amsthm}
\usepackage{tikz}
\usepackage{mathtools}
\usepackage{calc}
\usepackage{tikz}
\usepackage{pgfplots}
\usetikzlibrary{patterns,arrows,decorations.pathreplacing}
\usepackage{diagbox}
\usepackage{makecell}

\def\q{\hfill\rule{1ex}{1ex}}
\def\0{\emptyset}

\def\q{\hfill\rule{1ex}{1ex}}

\usepackage{amsthm}

\newtheorem{theorem}{Theorem}[section]

\newtheorem{lemma}[theorem]{Lemma}

\newtheorem{cor}[theorem]{Corollary}
\newtheorem{prop}[theorem]{Proposition}

\newtheorem{prob}[theorem]{Problem}

\renewcommand\theenumi{\rm \roman{enumi}}

\usepackage{algorithm}
\usepackage{algpseudocode}
\usepackage{graphicx}
\usepackage{amsmath}
\usepackage{pstricks,pgf,subcaption,amssymb}
\begin{document}
\title{\bf The Rainbow Saturation Number of Cycles
	 }
\author[]{
Yiduo Xu\thanks{E-mail:\texttt{xyd23@mails.tsinghua.edu.cn}}\;}
\author[2]{
Zhen He\thanks{Corresponding author. E-mail:\texttt{zhenhe@bjtu.edu.cn}}}
\author[]{
Mei Lu\thanks{E-mail:\texttt{lumei@mail.tsinghua.edu.cn}}}

\affil[]{\small Department of Mathematical Sciences, Tsinghua University, Beijing 100084, China}
\affil[2]{\small School of Mathematics and Statistics, Beijing Jiaotong University, Beijing 100044, China.}
\date{}

\maketitle\baselineskip 16.3pt

\begin{abstract}
	An edge-colored graph $(G,\mathcal C)$ is $F$-rainbow saturated if it contains no rainbow copy of $F$, but adding any nonedge in any prescribed color creates one.
	The rainbow saturation number $\operatorname{rsat}(n,F)$ is the minimum number of edges in such a graph of order $n$.
	In this paper we study $\operatorname{rsat}(n,C_r)$.
	We determine the value for $C_4$ exactly, and establish  bounds for every fixed $C_r$ with $r\geq5$.
\end {abstract}
	

\section{Introduction}

All graphs are finite, simple, and undirected.
We write $|G|$ and $e(G)$ for the numbers of vertices and edges of $G$, respectively, $d_G(v)$ for the degree of $v$, and $\delta(G)$ for the minimum degree.
We use $C_r$, $P_r$, and $K_r$ for the cycle, path, and complete graph on $r$ vertices, and $[k]=\{1,\ldots,k\}$.
All other graph-theoretic notation is standard.
	
Given graphs $G$ and $F$, we say  $G$ is \textit{$F$-saturated}, if $G$ does not contain any copy of $F$ but the addition of any nonedge $e \not \in E(G)$ would create a copy of $F$.
The \textit{saturation number} of $F$ is denoted by
\begin{equation*}
	   sat(n,F)=\min \{ e(G):G  \; \text{is}  \; F\text{-saturated  and } |G|=n  \} \, .
\end{equation*}
The first saturation problem was studied in 1964 by Erd\H os, Hajnal and Moon \cite{EHM}.
Since then, the topic has attracted extensive further research in the literature.
For readers interested in saturation problem, we refer to the survey \cite{survey}.

	An \textit{edge-colored graph} is a pair $(G,\mathcal C)$, where $\mathcal C:E(G)\to\mathbb N$.
	A subgraph $H\subseteq G$ is \textit{rainbow} if the restriction of $\mathcal C$ to $E(H)$ is injective.
	We call $(G,\mathcal C)$ \textit{$F$-rainbow free} if it contains no rainbow copy of $F$.
	The pair $(G,\mathcal C)$ is \textit{$F$-rainbow saturated} if it is $F$-rainbow free and adding any nonedge $uv$ in any color $\alpha\in\mathbb N$ creates a rainbow copy of $F$.
	Throughout, $\mathcal R$ denotes an injective edge-coloring. We call such a coloring a \textit{rainbow coloring}, and write $(G,\mathcal R)$ for the corresponding colored graph.
	The \textit{rainbow saturation number} of $F$ is
\begin{equation*}
	   \operatorname{rsat}(n,F)=\min \{ e(G): |G|=n
	   \text{ and $(G,\mathcal C)$ is $F$-rainbow saturated for some $\mathcal C$}\}.
\end{equation*}
	
	The rainbow saturation number was first studied by Gir\~ao, Lewis and Popielarz \cite{Gir}.
	Recently Behague, Johnston, Letzter, Morrison and Ogden \cite{BEH} determine that $rsat(n,F)=O(n)$ for any nonempty graph $F$, and gave an upper bound on $rsat(n,K_r)$.
	They also conjectured that there exists a constant $c=c(F)$ such that $rsat(n,F)=(c+o(1))n$ as $n \rightarrow \infty$.
	Chakraborti, Hendrey, Lund and Tompkins \cite{CHA} proved that there exists $\alpha_r$ such that $\lim \limits_{n \rightarrow \infty} \frac{rsat(n,K_r)}{n} = \alpha_r$ and showed the bounds of $\alpha_r$.

	In this paper we focus on rainbow saturation number of cycles.
	For saturation numbers of cycles there are a series of results.
	We refer the readers to \cite{FI,OLL,TUZ} for $sat(n,C_4)$, \cite{Chen1,Chen2,FI2} for $sat(n,C_5)$, \cite{Gou,LAN,Zhang} for $sat(n,C_6)$ and \cite{Fur,MPT} for $sat(n,C_r)$ where $r \geq 7$.
	An easy observation on a cycle rainbow saturated graph $G$ is that $G$ must be connected.
	Otherwise, the addition of two vertices in different components would not create a cycle, a contradiction.
	We first determined the exact value of $rsat(n,C_4)$.
\vspace{0.1em}	
	
\begin{theorem}\label{T11} For  $ n \geq 5$, $rsat(n,C_{4}) =3 \lceil \frac{n-1}{2} \rceil$.
\end{theorem}

	For $C_5$, we establish a lower bound by a weight-transfer argument, and an upper bound of leading coefficient $9/5$ by gluing sparse saturated blocks.
    % Even though the coloring on an optimal $C_5$-rainbow saturated graph $G$ is unknown, we are able to prove that if $G$ can be rainbow-colored then $e(G)$ would reach to the upper bound of $rsat(n,C_{5})$.
	
\vspace{0.1em}

\begin{theorem}\label{T12} For $n\geq6$,
\[
\frac32n-\frac52\leq \operatorname{rsat}(n,C_5)
\leq\left\lceil\frac{9(n-1)}5\right\rceil .
\]
\end{theorem}
\vspace{0.1em}

	For $C_6$, identifying copies of a fixed $13$-vertex block at a common cut vertex yields an upper bound of leading coefficient $11/6$, and for $r\geq7$ a periodic construction gives a linear upper bound with leading coefficient $5/3$.
\vspace{0.1em}

\begin{theorem}\label{T14} For $n\geq37$,
\[
\operatorname{rsat}(n,C_6)\leq 3(n-1)-14\left\lfloor\frac{n-1}{12}\right\rfloor = \frac{11}{6}n - O(1).
\]
\end{theorem}
\vspace{0.1em}

\begin{theorem}\label{T15}
For $r\geq7$ and $n\geq12(r-2)+1$,
\[
\operatorname{rsat}(n,C_r)
\leq \frac53(n-1)+4(r-2).
\]
\end{theorem}
\vspace{0.1em}

	Finally, we give a general lower bound for $\operatorname{rsat}(n,C_r)$ by characterizing the degree-two vertices of a $C_r$-rainbow saturated graph.
\vspace{0.5em}

\begin{theorem}\label{T17}
For $r\geq6$, define $F_0=1$, $F_1=3$, and $F_t=3+2F_{t-2}$ for $t\geq2$.
Then
\[
\operatorname{rsat}(n,C_r)
\geq\left(\frac65+\frac{1}{5F_{r-2}}\right)n-O(1).
\]
\end{theorem}
Here
\[
F_{2k}=2^{k+2}-3,\qquad F_{2k+1}=3\cdot2^{k+1}-3.
\]
%\vspace{0.0em}

	The rest of the paper is organized as follows.
	Section~2 records the necessary definitions and a color-avoidance criterion.
	In Section~3 we determine $\operatorname{rsat}(n,C_4)$.
	Sections~4 and~5 establish the bounds for $C_5$ and $C_6$, respectively.
	Section~6 gives a periodic construction for $C_r$ with $r\geq7$, and Section~7 proves a general lower bound.
	We collect open problems in Section~8.

\section{Preliminaries}

Let $(G,\mathcal C)$ be a colored graph.
	A subgraph $H$ of $G$ \textit{avoids} an edge $e$ if $e\notin E(H)$, and it avoids a color $\alpha$ if no edge of $H$ receives color $\alpha$.
For a nonedge $uv$, let $\mathcal F_{uv}$ be the family of copies $H$ of $F$ in $G+uv$ such that $uv\in E(H)$ and $H-uv$ is rainbow under $\mathcal C$.

\begin{prop}\label{P21}
Let $(G,\mathcal C)$ be $F$-rainbow free.
Then $(G,\mathcal C)$ is $F$-rainbow saturated if and only if, for every nonedge $uv$ and every color $\alpha\in\mathbb N$, there exists $H\in\mathcal F_{uv}$ such that $H-uv$ avoids $\alpha$.
\end{prop}

\noindent {\bf Proof. }If such an $H$ exists, then coloring the new edge $uv$ with $\alpha$ makes $H$ rainbow.
Conversely, suppose $(G,\mathcal C)$ is $F$-rainbow saturated.
After adding $uv$ in color $\alpha$, choose a rainbow copy $H$ of $F$ containing the new edge.
Then $H\in\mathcal F_{uv}$ and no edge of $H-uv$ has color $\alpha$. \qed
\vspace{0.0em}

	%We apply this criterion to cycles.
	
\begin{cor}\label{C22}(Sufficiency)
Let $(G,\mathcal C)$ be $C_r$-rainbow free.
Suppose that, for every nonedge $uv$, there are two edge-disjoint $u$--$v$ paths $P_r^1,P_r^2$ such that $P_r^1\cup P_r^2$ is rainbow.
Then $(G,\mathcal C)$ is $C_r$-rainbow saturated.
\end{cor}
\vspace{0.0em}

\begin{cor}\label{C23}(Necessity)
If $(G,\mathcal C)$ is $C_r$-rainbow saturated, then, for every nonedge $uv$ and every edge $e\in E(G)$, there exists a rainbow $u$--$v$ path of order $r$ that avoids $e$.
\end{cor}
%\vspace{0.0em}

\begin{cor}\label{C24}(Sufficiency and Necessity)
Let $\mathcal R$ be a rainbow coloring of $G$.
Then $(G,\mathcal R)$ is $C_r$-rainbow saturated if and only if:
\begin{enumerate}
\item $G$ is $C_r$-free;
\item for every nonedge $uv$ and every edge $e$, there exists a $u$--$v$ path of order $r$ avoiding $e$.
\end{enumerate}
\end{cor}
\vspace{0.0em}

\begin{lemma}\label{L25}
For $|G|\geq r\geq4$, if $(G,\mathcal C)$ is $C_r$-rainbow saturated, then $\delta(G)\geq2$.
\end{lemma}
\noindent {\bf Proof. }Suppose $\delta(G)=1$.
Let $u\in V(G)$ with $N(u)=\{v\}$ and let $w\in V(G)\setminus\{u,v\}$.
Every $u$--$w$ path of order $r$ contains $uv$, contradicting Corollary \ref{C23}. \qed

\section{Rainbow saturation of $C_4$ }

%\begin{figure}[H]
%\centering
%\begin{tikzpicture}[scale=.45]
%
%\draw (-10,0) arc(0:360:1cm and 3.5cm) ;
%\draw (-6,0) arc(0:360:1cm and 3.5cm) ;
%\draw (-2,0) arc(0:360:1cm and 3.5cm) ;
%
%\draw (-3,0) node[align=center]{$B_3$};
%\draw (-7,0) node[align=center]{$B_2$};
%\draw (-11,0) node[align=center]{$B_1$};
%
%\filldraw (-2.99,-5) circle (5pt);
%\filldraw (-6.99,-5) circle (5pt);
%\filldraw (-10.99,-5) circle (5pt);
%
%\draw (-3.2,-4.4) node[align=center]{$x_3$};
%\draw (-6.6,-4.4) node[align=center]{$x_2$};
%\draw (-10.6,-4.4) node[align=center]{$x_1$};
%
%\draw[dotted,thick] (-1, -5) arc(0:360:6cm and 1.2cm) ;
%%\draw (-7,-7) node[align=center]{$G[A]$};
%
%\draw (-2.99,-5)--(-6.99,-5) (-10.99,-5)--(-6.99,-5)  (-2.99,-5) arc(0:-180: 4cm and 0.7cm) (-10.99,-5);
%
%\draw (-10.5,-3) --(-6.99,-5) (-6.5,-3) --(-2.99,-5) (-3.5,-3) --(-10.99,-5);
%
%\end{tikzpicture}\\
%
%\caption{\centering $W^{(4,3)} = V_1 \cup V_2 \cup V_3$, $V_i = B_i \cup \{x_i \}$, $A=\{x_1,x_2,x_3\}$. The solid line represents the complete connection between vertices, and the dotted ellipse represents $W^{(4,3)}[A]= \gamma^{(3,3)}=K_3$. }
%
%\end{figure}

\begin{figure}[H]
%Figure 1
\centering
\begin{tikzpicture}[scale=2]

\draw[red,thick] (0,0) -- (0.5,1) ;
\draw[blue,thick] (0.5,1) -- (-0.5,1);
\draw[brown,thick] (-0.5,1) -- (0,0);

\draw[thick,green] (0,0) -- (0.5,-1) ;
\draw[thick,purple] (0.5,-1) -- (-0.5,-1) ;
\draw[thick,color=yellow!95!black!255!](-0.5,-1) -- (0,0);

\draw[thick,gray] (0,0) -- (1,0.5) ;
\draw[thick,cyan] (1,0.5) -- (1,-0.5);
\draw[thick,orange] (1,-0.5) -- (0,0);

\draw[thick,color=green!70!red!70!blue!70!] (0,0) -- (-1,0.5) ;
\draw[thick,color=green!80!red!20!blue!30!] (-1,0.5) -- (-1,-0.5) ;
\draw[thick,color=yellow!30!red!40!blue!40!] (-1,-0.5) -- (0,0);

\filldraw (0,0) circle (1pt);
\filldraw (0.5,1) circle (1pt);
\filldraw (-0.5,1) circle (1pt);
\filldraw (0.5,-1) circle (1pt);
\filldraw (-0.5,-1) circle (1pt);
\filldraw (1,0.5) circle (1pt);
\filldraw (1,-0.5) circle (1pt);
\filldraw (-1,0.5) circle (1pt);
\filldraw (-1,-0.5) circle (1pt);

\draw[red,thick] (4,0) -- (4.5,1) ;
\draw[blue,thick] (4.5,1) -- (3.5,1);
\draw[brown,thick] (3.5,1) -- (4,0);

\draw[red,thick] (4,0.66) -- (3.5,1) ;
\draw[blue,thick] (4,0) -- (4,0.66);
\draw[brown,thick] (4.5,1) -- (4,0.66);

\draw[thick,green] (4,0) -- (4.5,-1) ;
\draw[thick,purple] (4.5,-1) -- (3.5,-1) ;
\draw[thick,color=yellow!95!black!255!](3.5,-1) -- (4,0);

\draw[thick,gray] (4,0) -- (5,0.5) ;
\draw[thick,cyan] (5,0.5) -- (5,-0.5);
\draw[thick,orange] (5,-0.5) -- (4,0);

\draw[thick,color=green!70!red!70!blue!70!] (4,0) -- (3,0.5) ;
\draw[thick,color=green!80!red!20!blue!30!] (3,0.5) -- (3,-0.5) ;
\draw[thick,color=yellow!30!red!40!blue!40!] (3,-0.5) -- (4,0);

\filldraw (4,0.66) circle (1pt);
\filldraw (4,0) circle (1pt);
\filldraw (4.5,1) circle (1pt);
\filldraw (3.5,1) circle (1pt);
\filldraw (4.5,-1) circle (1pt);
\filldraw (3.5,-1) circle (1pt);
\filldraw (5,0.5) circle (1pt);
\filldraw (5,-0.5) circle (1pt);
\filldraw (3,0.5) circle (1pt);
\filldraw (3,-0.5) circle (1pt);

\draw (0,-0.2) node[align=center]{$u$};
\draw (-0.7,1) node[align=center]{$x_1$};
\draw (0.7,1) node[align=center]{$x_2$};
\draw (-1.2,0.5) node[align=center]{$y_1$};
\draw (-1.2,-0.5) node[align=center]{$y_2$};

\draw (4,-0.2) node[align=center]{$u$};
\draw (4,0.8) node[align=center]{$x_3$};
\draw (3.3,1) node[align=center]{$x_1$};
\draw (4.7,1) node[align=center]{$x_2$};
\draw (2.8,0.5) node[align=center]{$y_1$};
\draw (2.8,-0.5) node[align=center]{$y_2$};

\end{tikzpicture}\\

\caption{\centering The coloring of $M_n$.}

\end{figure}

For $n\geq5$, let
\[
M_n=
\begin{cases}
K_1\vee\bigl(\frac{n-1}{2}K_2\bigr) & \text{if $n$ is odd},\\[0.3em]
K_1\vee\bigl(K_3\cup\frac{n-4}{2}K_2\bigr) & \text{if $n$ is even}.
\end{cases}
\]
When $n$ is odd we give $M_n$ a rainbow coloring $\mathcal C=\mathcal R$.
When $n$ is even we start from a rainbow coloring and recolor the edges of the unique $K_4$ by a proper $3$-coloring (see Figure~1).
We will show that $(M_n,\mathcal C)$ is $C_4$-rainbow saturated.

	Clearly $(M_n,\mathcal C)$ is $C_4$-rainbow free.
	By symmetry we only need to consider the addition of nonedge $x_1y_1$ (Figure 1).
	Since $x_1uy_2y_1$ and $x_1x_2uy_1$ are two edge-disjoint $P_4$'s connecting $x_1,y_1$ whose union is rainbow, by Corollary \ref{C22}, we know that $(M_n,\mathcal C)$ is $C_4$-rainbow saturated.
	Thus $rsat(n,C_4) \leq e(M_n) = 3 \lceil \frac{n-1}{2} \rceil$.
	The following proposition shows that the lower bound of $rsat(n,C_4)$ is also $3 \lceil \frac{n-1}{2} \rceil$ and we finish the proof of Theorem \ref{T11}.
\vspace{0.4em}	
	
\begin{prop}\label{P31}
For $n \geq 5$, $rsat(n,C_4) \geq 3 \lceil \frac{n-1}{2} \rceil$.
\end{prop}

\noindent {\bf Proof. }Let $(G,\mathcal C)$ be a $C_4$-rainbow saturated graph of order $n$.
If $\delta(G) \geq 3$ then we are done.
Assume that $u \in V(G)$ is a degree 2 vertex.
Let $N(u)=\{u_1,u_2\}$, $A=(N(u_1) \cup N(u_2)) \setminus \{u,u_1,u_2\}$ and $B=V(G) \setminus (A \cup \{u,u_1,u_2\})$.
Then $|A|+|B|=n-3$ and $N(v) \cap A\not=\emptyset$ for any $v\in A$ by Corollary \ref{C23} (otherwise all $P_4$'s connecting $v,u$ can not avoid $u_1u_2$).

	For any $v \in B$, we have $|N(v) \cap A| \geq 2$.
	Otherwise all $P_4$'s connecting $u,v$ can not avoid $vw$, where $w \in N(v) \cap A$, a contradiction to Corollary \ref{C23}.
	We consider the following two cases.
\vspace{0.4em}	
	
\noindent{\bf Case 1. }$u_1u_2 \not \in E(G)$.
Then for any $v \in A$, $|N(v) \cap A| \geq 2$.
Otherwise all $P_4$'s connecting $u,v$ can not avoid $vw$, where $w \in N(v) \cap A$, a  contradiction to Corollary \ref{C23}.
So $e(A) \geq |A|$.
Thus
\begin{equation}
\begin{aligned}
e(G) & = e(G[ \{ u,u_1,u_2 \} ]) + e(N(u),A) + e(A) + e(A,B) +e(B) \\
& \geq 2 + |A| + |A|+2|B| + 0 = 2n-4.
\end{aligned}
\end{equation}
For $n \geq 5$ and $n \neq 6$ we have $2n-4 \geq 3 \lceil \frac{n-1}{2} \rceil$.
For $n=6$, we have $|N(u_i)\cap A|\ge 2$ for $i\in [2]$ by Corollary \ref{C23} and $u_1u_2\notin E(G)$.
Hence $e(G)\ge 9= 3 \lceil \frac{n-1}{2} \rceil.$
%an easy argument can show that the equality can not hold in (1) and we left the details to readers.
Therefore in this case we have $e(G) \geq 3 \lceil \frac{n-1}{2} \rceil$.
\vspace{0.4em}

\noindent{\bf Case 2. }$u_1u_2  \in E(G)$.
Let $v \in A$ and assume that $u_1 \in N(v)$.
Since there must exist a $P_4$ connecting $u,v$ avoiding $vu_1$ by Corollary \ref{C23}, we have $u_2 \in N(v)$ or $N(v) \cap A \neq \emptyset$.
Thus we have $e(N(u),A)+e(A) \geq |A|+\frac{1}{2}|A|$ and then
\begin{equation}
\begin{aligned}
e(G) & = e(G[ \{ u,u_1,u_2 \} ]) + e(N(u),A) + e(A) + e(A,B) +e(B) \\
& \geq 3 + |A| + \frac{1}{2}|A|+2|B| + 0 \geq  \frac{3n-3}{2}.
\end{aligned}
\end{equation}
We are done if $n$ is odd.
If $n$ is even, we have $e(G) \geq \frac{3}{2}n-1$.
Suppose $e(G) = \frac{3}{2}n-1$.
Then $|B| \leq 1$ and $e(A,B)=2|B|$.
Moreover if $|B|=1$ then $e(N(u),A)=|A|, \, e(A)=\frac{1}{2}|A|$, else if $|B|=0$ then either $e(N(u),A)=|A|+1, \, e(A)=\frac{1}{2}(|A|-1) $ or $e(N(u),A)=|A|, \, e(A)=\frac{1}{2}(|A|+1) $.
Hence there exists  $w \in V(G)\setminus \{u,u_1,u_2\}$ such that

	(i) 	$A' \subseteq A$, $G[A']$ is a perfect matching and any vertex $v\in A'$ is located in a triangle together with some $v' \in A' , u_i \in N(u)$ (by Corollary \ref{C23}), where $A'=V(G) \setminus \{ u , u_1,u_2,w\}$;
	
	(ii)\, either $w \in B$ with $N(w)=\{v_1,v_2\}$ for some $v_1,v_2 \in A'$, or $w \in A$ with $N(w) = \{u_1,u_2\}$, or $w \in A$ with $N(w)=\{u_i, v_1 \}$ for some $v_1 \in A'$ and $i \in [2]$.
\vspace{0.4em}

%%%%%%%%%%%%%%%%%%%%%%%%%%%%%%%%%%%%%%%
\noindent{\bf Claim 1. }$B =\emptyset$, i.e. $w \not \in B$.

\noindent{\bf Proof of Claim 1.} Suppose that $w \in B$ and $N(w)=\{v_1,v_2\} \subseteq A'=A$.
Since there exists a $P_4$ connecting $w,u$ avoiding $uu_i$ for $i \in [2]$, $N(v_1) \cap N(u) \neq N(v_2) \cap N(u)$.
Assume that $u_iv_i \in E(G)$, $i=1,2$.
By (i), $v_1v_2\notin E(G)$.
Then there exists only one $P_4=v_1u_1u_2v_2$ connecting $v_1,v_2$, a contradiction to Corollary \ref{C23}. \q
\vspace{0.4em}

\noindent{\bf Claim 2. }There exists $i \in [2]$ such that for all $v \in A'$, $N(v) \cap N(u) =\{u_i\}$.

\noindent{\bf Proof of Claim 2.}	Suppose  that there exist $v_1,v_2 \in A'$ such that $v_iu_i \in E(G)$, $i \in [2]$.
By (i), $v_1v_2\notin E(G)$.
By (ii) and Claim 1, we can assume $wv_2\notin E(G)$.
Then there is no  $P_4$ connecting $v_1,v_2$  avoiding   $v_2u_2$, a contradiction to Corollary \ref{C23}. \q
\vspace{0.4em}
	
	By Claim 2, we can assume that all vertices in $A'$ are connected to $u_1$.
	By Claim 1, $w\in A$.
	If  $N(w) = \{u_1,u_2\}$, then all $P_4$'s connecting $u,w$ can not avoid $u_1u_2$, a contradiction to Corollary \ref{C23}.
	Suppose  $N(w)=\{u_i, v_1 \}$ for some $v_1 \in A'$ and $i \in [2]$.
	Let $v_1'\in N(v_1)\cap A'$.
	Then for $i=1$ there exists no $P_4$ connecting $w,v_1'$ avoiding $u_1v_1$, and for $i=2$ there exists no $P_4$ connecting $u,w$ avoiding $uu_1$, a contradiction to Corollary \ref{C23}.
	Hence $rsat(n,C_4) \geq 3 \lceil \frac{n-1}{2} \rceil$ holds when $n$ is even.  \qed

\section{Rainbow saturation of $C_5$ }	

%	We first construct a graph which is $C_5$-rainbow saturated under a rainbow coloring.

	For $t\geq3$, let $W_t=K_2\vee\overline{K_{t-2}}$ and write $u_1,u_2$ for the two vertices of degree $t-1$, then $W_t$ is $C_5$-free.
	Give $W_t$ a rainbow coloring $\mathcal R$.
	If $t\geq6$, let $I=\{x_1,\ldots,x_{t-2}\}$ be its independent set.
	For distinct $x_i,x_j\in I$, choose distinct $k,\ell\in[t-2]\setminus\{i,j\}$.
	Then
	\(
	x_iu_1x_ku_2x_j
	\text{ and }
	x_iu_2x_\ell u_1x_j
	\)
	are edge-disjoint copies of $P_5$.
	Hence Corollary~\ref{C22} shows that $(W_t,\mathcal R)$ is $C_5$-rainbow saturated for every $t\geq6$.
	Propositions~\ref{P5new} and \ref{P41} below prove Theorem~\ref{T12}.
	
\begin{figure}[H]
%Figure 2
\centering
\begin{tikzpicture}[scale=2]

\draw[dotted,thick] (5,0) arc(0:360: 1.7cm and 0.5cm) ;

\draw[thick,color=green!80!red!20!blue!30!] (1,0.8) -- (1,-0.8);

\draw[red,thick] (1,0.8) -- (2,0) ;
\draw[blue,thick] (1,0.8) -- (2.5,0);
\draw[brown,thick] (1,0.8) -- (4.5,0);

\draw[thick,green] (1,-0.8) -- (2,0) ;
\draw[thick,purple] (1,-0.8) -- (2.5,0) ;
\draw[thick,cyan] (1,-0.8) -- (4.5,0);

\filldraw (1,0.8) circle (1pt);
\filldraw (1,-0.8) circle (1pt);
\filldraw (2,0) circle (1pt);
\filldraw (2.5,0) circle (1pt);
\filldraw (3,0) circle (0.5pt) (3.2,0) circle (0.5pt) (3.4,0) circle (0.5pt)  (3.6,0) circle (0.5pt) (3.8,0) circle (0.5pt) (4,0) circle (0.5pt);
\filldraw (4.5,0) circle (1pt);

\draw (3.3,-0.7) node[align=center]{$I$};
\draw (1,1) node[align=center]{$u_1$};
\draw (1,-1) node[align=center]{$u_2$};

\end{tikzpicture}\\

\caption{\centering The rainbow-colored graph $(W_t,\mathcal R)$.}

\end{figure}

\begin{prop}\label{P5new}
For every $n\geq6$,
\[
\operatorname{rsat}(n,C_5)
\leq\left\lceil\frac{9(n-1)}5\right\rceil .
\]
\end{prop}

\noindent {\bf Proof. }Put $q=\lfloor(n-1)/5\rfloor$ and write $n-1=5q+s$, where $0\leq s\leq4$.
Take $q-1$ copies of $W_6$ and one copy of $W_{6+s}$.
In each copy choose a vertex from the independent set and identify all the chosen vertices to a single vertex $x$.
Denote the resulting graph by $G_n$.
Thus distinct blocks meet only at the cut vertex $x$.
Give all edges of $G_n$ distinct colors.

Every cycle of $G_n$ is contained in one block, and each block is a copy of some $W_t$, hence $C_5$-free.
Thus $G_n$ is $C_5$-free.
A nonedge whose endpoints lie in one block is handled inside that block by the construction of $W_t$.

It remains to consider a nonedge with endpoints in two distinct blocks.
In each block $i$, let $a_i,b_i$ be the two vertices with maximum degree inside the block, and let $p_i,q_i \neq x$ be two vertices of the independent set of the block.
Up to symmetry it suffices to check the three endpoint types below. In each type no edge of $G_n$ lies in every listed $P_5$\,:
\[
\begin{array}{|c|l|}
\hline
\text{endpoints}&\text{$P_5$'s}\\ \hline
p_1,p_2&
p_1a_1xa_2p_2,\quad p_1b_1xb_2p_2\\ \hline
p_1,a_2&
p_1a_1xb_2a_2,\quad
p_1a_1b_1xa_2,\quad
p_1b_1xb_2a_2\\ \hline
a_1,a_2&
a_1xb_2q_2a_2,\quad
a_1q_1b_1xa_2 .\\ \hline
\end{array}
\]
All these paths are rainbow under $\mathcal R$.
%Consequently, for every nonedge and every edge of $G_n$, one of the listed paths avoids that edge.
Corollary~\ref{C24} shows that $(G_n,\mathcal R)$ is $C_5$-rainbow saturated.

Finally,
\[
e(G_n)=\sum_{i=1}^q(2t_i-3)
=2\sum_{i=1}^q(t_i-1)-q
=2n-2-q.
\]
If $n-1=5q+s$, then $2n-2-q=9q+2s=\lceil9(n-1)/5\rceil$, as required. \qed
\vspace{0.5em}

\begin{prop}\label{P41}
For $n \geq 5$, $rsat(n,C_5) \geq \frac{3}{2}n-\frac{5}{2}$.
\end{prop}

\noindent {\bf Proof. }Let $(G,\mathcal C)$ be a $C_5$-rainbow saturated graph of order $n$.
If $\delta(G) \geq 3$ then we immediately have the result.
So we assume that there exists a degree 2 vertex $u \in V(G)$.
Since $(G,\mathcal C)$ is $C_5$-rainbow saturated, $diam(G) \leq 4$.
Let $S_0=\{u\}$ and $S_i=\{v:d(u,v)=i\}$ for $i \in [4]$.
Then $S_1=\{u_1,u_2\}$.
For $i \in [4]$, let $f$ be a weight function such that
\begin{equation*}
 f(v)=|N_{i-1}(v)|+\frac{1}{2}|N_i(v)| \, , \quad \forall \, v \in S_i \, ,
\end{equation*}
where $N_j(v)=N(v) \cap S_j$.
Clearly we have $\sum \limits_{v \in V(G) \setminus \{u\}} f(v) = e(G)$ since all edges of $G$ have been calculated for exactly one time.

	For every $v\in S_i$ with $i\geq1$, connectivity gives $|N_{i-1}(v)|\geq1$.
	Moreover: (1) $f(v)=1$ if and only if $|N_{i-1}(v)|=1$ and $|N_i(v)|=0$; (2) $f(v)=\frac{3}{2} $ iff $|N_{i-1}(v)|=1$ and  $|N_{i}(v)|=1$ ; (3) $f(v) \geq 2$ iff $|N_{i-1}(v)| \geq 2$ or  $|N_{i}(v)|\geq 2$.
	Let $A_2 = \{v : v \in S_2, |N_1(v)|=1, |N_2(v)|=0\}$ and $B_2 = S_2 \setminus A_2$.
	Let $A_3 = \{v : v \in S_3 , N_2(v) \cap A_2 \neq \emptyset \text{ and } |N_2(v)| \geq 2\}$ and $B_3 = S_3 \setminus A_3$.
	Clearly for $v \in B_2 $ we have $f(v) \geq \frac{3}{2}$ and for $v \in  S_1$ we have $f(v) \geq 1$.
	Also for $v\in A_3$, $f(v)\ge 2$.
\vspace{0.5em}

\noindent {\bf Claim 1.} For every $v \in S_4$, $f(v) \geq |N_{3}(v)| \geq 2$.

\noindent {\bf Proof of Claim 1.} Suppose there is $v \in S_4$ such that $N_3(v)=\{w\}$.
Since $d(u,v)=4$, every $u$--$v$ path of order $5$ contains $vw$, contradicting Corollary~\ref{C23}. \q
\vspace{0.5em}

\noindent {\bf Claim 2.} For every $v \in B_3$, $f(v) \geq \frac{3}{2}$.

\noindent {\bf Proof of Claim 2.} Suppose there is $v\in B_3$ with $f(v)=1$.
Then $N_3(v)=\emptyset$ and $N_2(v)=\{w\}$ for some $w$.
In every copy $P_5=ux_1x_2x_3v$, the vertex $x_3$ does not lie in $S_4$, and hence $x_3=w$.
Thus every such path contains $vw$, contradicting Corollary~\ref{C23}. \q
\vspace{0.5em}

\noindent {\bf Claim 3.} For every $v \in A_2$, $N_3(v) \cap A_3 \neq \emptyset$.

\noindent {\bf Proof of Claim 3.}	Let $v \in A_2$.
Then $|N_1(v)|=1$ and $N_2(v)=\emptyset$.
Assume $N_1(v)=\{u_1\}$.
Since there exists a $P_5=ux_1x_2x_3v$ connecting $u,v$ avoiding $vu_1$ and $N_2(v)=\emptyset$, we must have $x_i \in S_i$ for $i \in [3]$.
Since $v,x_2\in S_2$, $x_3 \in N_3(v) \cap A_3$.  \q
\vspace{0.5em}

	 For $v \in A_3$, let $q(v) = |N_2(v) \cap A_2|$.
	 We partition $A_3$ into $A_3^+=\{ v \in A_3 : f(v) \geq \frac{5}{2} \}$ and $A_3^-=A_3 \setminus A_3^+$.
	 Denote $f_1(v)=f(v)-\frac{1}{2}q(v)$ for $v \in A_3^+$, $f_1(v)=f(v)$ for $v \in A_3^-$ and $f_1(w)=f(w)+\frac{1}{2}|N(w) \cap A_3^+|$ for $w \in A_2$.
	 Clearly $\sum_{v \in A_2 \cup A_3} f(v) = \sum_{v \in A_2 \cup A_3} f_1(v)$.
	 Partition $A_2$ into $A_2^+=\{ w \in A_2 : f_1(w) > f(w) \}$ and $A_2^-=A_2 \setminus A_2^+$.
	 Still it is obvious that for any $w \in A_2^+$, we have $f_1(w) \geq \frac{3}{2}$.
\vspace{0.5em}

\noindent {\bf Claim 4.} For any $v \in A_3^+$, $f_1(v) \geq \frac{3}{2}$.

\noindent {\bf Proof of Claim 4.} Let $v \in A_3^+$.
Then $q(v) \geq 1$, $f(v)=q(v)+|N(v) \cap B_2|+\frac{1}{2}|N_3(v)|$ and $f(v) \geq \frac{5}{2}$.
So $f_1(v)= \frac{q(v)}{2}+|N(v) \cap B_2|+\frac{1}{2}|N_3(v)|$.
If $|N(v) \cap B_2| \geq 1$, then $f_1(v) \geq  \frac{3}{2}$ and we are done.
Assume $|N(v) \cap B_2|= 0$.
Then $f(v)=q(v)+\frac{1}{2}|N_3(v)|$ and $f_1(v)= \frac{q(v)}{2}+\frac{1}{2}|N_3(v)|$.
If $|N_3(v)| \geq 2$, then $f_1(v) \geq  \frac{3}{2}$ and we are done.
If $|N_3(v)| = 1$, then $q(v) \geq 2$ by $f(v) \geq \frac{5}{2}$ which implies $f_1(v) \geq  \frac{3}{2}$.
If $|N_3(v)| = 0$, then $q(v) \geq 3$ by $f(v) \geq \frac{5}{2}$ which implies $f_1(v) \geq  \frac{3}{2}$. \q
\vspace{0.5em}

	Next we focus on $A_3^-$ and $A_2^-$.
	Note that $f(v)=|N_2(v)|+ \frac{1}{2}|N_3(v)|$ and $|N_2(v)|\ge 2$ for $v\in A_3$.
	Then for $v \in A_3^-$, we have $|N_2(v)|=2$ and $|N_3(v)|=0$ which implies $N(v) \subseteq S_4 \cup S_2$.
	We assert that it is possible to send a weight $\frac{1}{2}$ from $A_3^-$ to $A_2^-$.
	By Claim 3 and the definition of $A_2^-$, for every vertex $w \in A_2^-$, $\emptyset\not=N(w)\cap A_3\subseteq A_3^-$.
\vspace{0.5em}
	
\noindent {\bf Claim 5.} For any $v \in A_3^-$, $|N(v) \cap A_2^-| \leq 1$.

\noindent {\bf Proof of Claim 5.}	Suppose for contradiction, let $\{w_1,w_2\} = N_2(v) \cap A_2^-=N_2(v)$.
Pick a $P_5$ connecting $u,v$ arbitrarily, says $vx_1x_2x_3u$.
Since $N(v) \subseteq S_4 \cup S_2$, we have  $x_1 \in \{w_1,w_2\}$.
Since $w_i\in A_2^-\subseteq A_2$,  $N(w_i) \in S_1 \cup S_3$ for $i \in [2]$.
So we have $x_2 \in S_1$.
Moreover $x_3 \in S_1$ and hence $x_2x_3 = u_1u_2$.
Thus all $P_5$'s connecting $u,v$ can not avoid $u_1u_2$, a contradiction to Corollary \ref{C23}. \q
\vspace{0.5em}	
	
	Denote $f_2(v)=f(v)-\frac{1}{2}|N(v) \cap A_2^-|$ for $v \in A_3^-$ and $f_2(w)=f(w)+\frac{1}{2}|N(w) \cap A_3^-|$ for $w \in A_2^-$.
	By Claim 5 and $f(v)\ge 2$ for $v \in A_3^-$, we immediately have $f_2(v) \geq \frac{3}{2}$ for $v \in A_3^-$ and also, $f_2(w) \geq \frac{3}{2}$ for $w \in A_2^-$.
	Still obviously $\sum_{v \in A_2^- \cup A_3^- } f(v) = \sum_{v \in A_2^-  \cup A_3^- } f_2(v)$.
	
	 Recall that $|S_0|=1, |S_1|=2$ and $|A_2|+|B_2|+|A_3|+|B_3|+|S_4|=n-3$.
	 Then
\begin{equation*}
\begin{aligned}
	e(G)  & =  \sum_{v \in V(G) \setminus \{u\}} f(v) \\
	& =\sum_{v \in S_1 \cup B_2 \cup B_3 \cup S_4 } f(v) + \sum_{v \in A_2^- \cup A_3^-} f_2(v) + \sum_{v \in A_2^+ \cup A_3^+} f_1(v) \\
	& \geq |S_1|+ \frac{3}{2}(|B_2|+|B_3|)+ 2|S_4|+ \frac{3}{2} \left( |A_2| + |A_3| \right) \\
	& \geq 2 + \frac{3}{2}(n-3) = \frac{3}{2}n-\frac{5}{2}.
\end{aligned}
\end{equation*}
The proof is done. \qed

\section{Rainbow saturation of $C_6$}

\begin{figure}[H]
\centering
\begin{tikzpicture}[scale=1.5]

% Central K_5: h at the centre, u_1..u_4 around a square
% Spokes
\draw[thick,red] (0,0) -- (0,1.1);
\draw[thick,blue] (0,0) -- (1.1,0);
\draw[thick,green!70!black] (0,0) -- (0,-1.1);
\draw[thick,orange] (0,0) -- (-1.1,0);
% Outer 4-cycle
\draw[thick,purple] (0,1.1) -- (1.1,0);
\draw[thick,cyan] (1.1,0) -- (0,-1.1);
\draw[thick,brown] (0,-1.1) -- (-1.1,0);
\draw[thick,gray] (-1.1,0) -- (0,1.1);
% Diagonals (bent so as not to pass through h)
\draw[thick,color=yellow!80!red] (0,1.1) to[bend left=28] (0,-1.1);
\draw[thick,color=green!50!blue] (1.1,0) to[bend left=28] (-1.1,0);

% Petal at u_1 (top)
\draw[thick,color=red!70!black] (0,1.1) -- (-0.4,1.85);
\draw[thick,color=blue!70!black] (0,1.1) -- (0.4,1.85);
\draw[thick,color=green!60!black] (-0.4,1.85) -- (0.4,1.85);

% Petal at u_2 (right)
\draw[thick,color=orange!80!black] (1.1,0) -- (1.85,-0.4);
\draw[thick,color=purple!80!black] (1.1,0) -- (1.85,0.4);
\draw[thick,color=cyan!80!black] (1.85,-0.4) -- (1.85,0.4);

% Petal at u_3 (bottom)
\draw[thick,color=brown!80!black] (0,-1.1) -- (-0.4,-1.85);
\draw[thick,color=gray!80!black] (0,-1.1) -- (0.4,-1.85);
\draw[thick,color=yellow!70!red!60!black] (-0.4,-1.85) -- (0.4,-1.85);

% Petal at u_4 (left)
\draw[thick,color=green!40!blue!60] (-1.1,0) -- (-1.85,-0.4);
\draw[thick,color=red!50!blue] (-1.1,0) -- (-1.85,0.4);
\draw[thick,color=yellow!60!blue] (-1.85,-0.4) -- (-1.85,0.4);

% Vertices
\filldraw (0,0) circle (1.2pt);
\filldraw (0,1.1) circle (1.2pt);
\filldraw (1.1,0) circle (1.2pt);
\filldraw (0,-1.1) circle (1.2pt);
\filldraw (-1.1,0) circle (1.2pt);
\filldraw (-0.4,1.85) circle (1.2pt);
\filldraw (0.4,1.85) circle (1.2pt);
\filldraw (1.85,-0.4) circle (1.2pt);
\filldraw (1.85,0.4) circle (1.2pt);
\filldraw (-0.4,-1.85) circle (1.2pt);
\filldraw (0.4,-1.85) circle (1.2pt);
\filldraw (-1.85,-0.4) circle (1.2pt);
\filldraw (-1.85,0.4) circle (1.2pt);

% Labels
\draw (0.15,-0.15) node{$h$};
\draw (0.12,1.1) node[right]{$u_1$};
\draw (1.1,0.12) node[above]{$u_2$};
\draw (0.12,-1.1) node[right]{$u_3$};
\draw (-1.1,0.12) node[above]{$u_4$};
\draw (-0.4,2.05) node{$x_1$};
\draw (0.4,2.05) node{$y_1$};
\draw (2.1,-0.4) node{$x_2$};
\draw (2.1,0.4) node{$y_2$};
\draw (-0.4,-2.05) node{$x_3$};
\draw (0.4,-2.05) node{$y_3$};
\draw (-2.1,-0.4) node{$x_4$};
\draw (-2.1,0.4) node{$y_4$};

\end{tikzpicture}

\caption{\centering The rainbow-colored graph $(B,\mathcal R)$.}
\label{fig:B}
\end{figure}

Let $B$ be obtained from a copy of $K_5$ on $\{h,u_1,u_2,u_3,u_4\}$ by adding, for every $i\in[4]$, two vertices $x_i,y_i$ and the triangle $u_ix_iy_i$ (Figure~\ref{fig:B}).
Then $B$ has $13$ vertices and $22$ edges.

For $q\geq1$, take $q$ copies of $B$ and identify their vertices labelled $h$.
Call the resulting graph $B_q$, and give all its edges distinct colors.
Every cycle is contained in one block, and every cycle in a block is either contained in the central $K_5$ or in one of the attached triangles.
In particular, $B_q$ is $C_6$-free.
\vspace{0.2em}

\begin{lemma}\label{L6block}
The rainbow-colored graph $(B_q,\mathcal R)$ is $C_6$-rainbow saturated.
\end{lemma}
\noindent {\bf Proof. }By Corollary~\ref{C24}, it is enough to show that, for every nonedge $vw$ and every edge $e$, there is a $v$--$w$ path of order $6$ avoiding $e$.
By symmetry we only need to consider the following six types.
Primed vertices belong to a second block.
In each type no edge of $B$ occurs in all the displayed paths.
\[
\begin{array}{|c|l|}
\hline
\text{endpoints}&\text{$P_6$'s}\\ \hline
h,x_1&
hu_4u_2u_1y_1x_1,\quad hu_2u_3u_4u_1x_1\\ \hline
x_1,u_2&
x_1u_1u_4hu_3u_2,\quad x_1y_1u_1u_3u_4u_2\\ \hline
x_1,x_2&
x_1u_1u_4u_2y_2x_2,\quad x_1y_1u_1hu_2x_2\\ \hline
u_1,u'_1&
u_1u_4u_2hu'_2u'_1,\quad u_1u_2u_3u_4hu'_1\\ \hline
u_1,x'_1&
u_1u_4u_2hu'_1x'_1,\quad u_1hu'_2u'_1y'_1x'_1\\ \hline
x_1,x'_1&
x_1u_1u_4hu'_1x'_1,\quad
x_1u_1hu'_2u'_1x'_1,\quad
x_1u_1hu'_1y'_1x'_1,\quad
x_1y_1u_1hu'_1x'_1 .\\ \hline
\end{array}
\]
All displayed paths are rainbow under $\mathcal R$.
Thus at least one of them avoids the prescribed edge $e$. 
Hence $(B_q,\mathcal R)$ is $C_6$-rainbow saturated. \qed
\vspace{0.5em}

\noindent {\bf Proof of Theorem~\ref{T14}. }Write $n-1=12q+s$, where $0\leq s\leq11$.
Since $n\geq37$, we have $q\geq3$, and $B_q$ contains at least $s$ attached triangles.
Choose $s$ of them and replace each triangle on $\{u_i,x_i,y_i\}$ by a $K_4$ on $u_i$ and three petal vertices; each replacement adds one vertex and three edges.
Denote the resulting graph by $B_q^s$.
By an argument similar to that of Lemma~\ref{L6block}, $(B_q^s,\mathcal R)$ is $C_6$-rainbow saturated.
It has $n$ vertices and
\[
e(B_q^s)=22q+3s=3(n-1)-14\left\lfloor\frac{n-1}{12}\right\rfloor = \frac{11}{6}n-O(1).
\]\qed

\section{Rainbow saturation of $C_r$ for $r \geq 7$ }	
Put $d=r-2$.
Let $p\geq4d$ be divisible by $d$, and let $v_0v_1\cdots v_{p-1}v_0$ be a cycle (the \emph{rim}).
Indices are taken modulo $p$.
Add a vertex $h$ adjacent to every $v_i$ and call each edge $hv_i$ a \emph{spoke}.
For every $i$, add two vertices $x_i,y_i$ and the triangle $v_ix_iy_i$, and call each edge of this triangle a \emph{petal edge}.
Denote the resulting graph by $Q_{r,p}$ (Figure~\ref{fig:Qrp}).
Color the spokes periodically by $\mathcal C(hv_i)=i\pmod d$, and give every rim edge and every petal edge a distinct new color, disjoint from the spoke colors.

\begin{figure}[H]
\centering
\begin{tikzpicture}[scale=1.8]
% Schematic with d=5, p=20=4d.
% Spoke colour of hv_i is i mod 5; all other edges black.

\coordinate (h) at (0,0);

% All rim vertices v_0,...,v_19 (angle = 90 - 18k)
\foreach \k in {0,...,19} {
  \coordinate (v\k) at ({90-18*\k}:2.2);
}

% Black rim cycle
\draw[thick,black]
  (v0)--(v1)--(v2)--(v3)--(v4)--(v5)--(v6)--(v7)--(v8)--(v9)--(v10)
  --(v11)--(v12)--(v13)--(v14)--(v15)--(v16)--(v17)--(v18)--(v19)--(v0);

% Periodic spokes
\foreach \k/\col in {
  0/red,5/red,10/red,15/red,
  1/blue,6/blue,11/blue,16/blue,
  2/{green!70!black},7/{green!70!black},12/{green!70!black},17/{green!70!black},
  3/orange,8/orange,13/orange,18/orange,
  4/purple,9/purple,14/purple,19/purple} {
  \draw[thick,\col] (h)--(v\k);
}

% Slim black petal triangles at every rim vertex
\foreach \k in {0,...,19} {
  \pgfmathsetmacro{\ang}{90-18*\k}
  \coordinate (x\k) at ({\ang+5.5}:2.78);
  \coordinate (y\k) at ({\ang-5.5}:2.78);
  \draw[thick,black] (v\k)--(x\k)--(y\k)--(v\k);
  \filldraw (x\k) circle (0.7pt);
  \filldraw (y\k) circle (0.7pt);
}

% Vertices
\filldraw (h) circle (1.5pt);
\foreach \k in {0,...,19} {
  \filldraw (v\k) circle (1.1pt);
}

% Labels of v_i: inside the circle, clockwise of the spoke hv_i
\draw (0.28,-0.18) node{$h$};
\foreach \k in {0,...,19} {
  \pgfmathsetmacro{\ang}{90-18*\k-4}
  \draw (\ang:1.99) node{\scriptsize $v_{\k}$};
}
\draw (95.5:2.95) node{\scriptsize $x_0$};
\draw (84.5:2.95) node{\scriptsize $y_0$};
\draw (77.5:2.95) node{\scriptsize $x_1$};
\draw (66.5:2.95) node{\scriptsize $y_1$};

\end{tikzpicture}

\caption{\centering The graph $Q_{r,p}$, illustrated with $d=5$ and $p=20=4d$.
Spokes of equal residue modulo $d$ receive the same colour (e.g.\ $hv_0$, $hv_5$, $hv_{10}$ and $hv_{15}$);
every rim edge and every petal edge is rainbow and is drawn in black.}
\label{fig:Qrp}
\end{figure}

\begin{lemma}\label{Lperiodic}
For every nonedge $zw$ of $Q_{r,p}$ and every color $\alpha$, there is a rainbow $z$--$w$ path of order $r$ that avoids $\alpha$.
\end{lemma}

\noindent {\bf Proof. }Indices are taken modulo $p$.
By symmetry we only need to consider the following four types of nonedge
\(
	zw=hx_0, v_0v_t,v_0x_t, x_0x_t
\)
where $2\le t\le p/2$ if $zw=v_0v_t$ and $1\le t\le p/2$ if $zw=v_0x_t,x_0x_t$.

Clearly we only need to consider avoiding a color that appears in \(Q_{r,p}\).
If \(\alpha\) is used in \(Q_{r,p}\), then it is the color of a spoke, a rim edge, or a petal edge.
A spoke color is the common color of all spokes $hv_{i+md}$ for $0 \le m \le 3$, and is avoided precisely by avoiding all these spokes.

For $a,b\ge0$, $\sigma\in\{+1,-1\}$ and a rim index $t$, write $P^t_\sigma(a,b)$ for the $v_0$--$v_t$ path
\[
v_0v_{\sigma}v_{2\sigma}\cdots v_{a\sigma}\,h\,v_{t-\sigma b}v_{t-\sigma(b-1)}\cdots v_t
\]
If $a=0$ this is $v_0hv_{t-\sigma b}\cdots v_t$; if $b=0$ this is $v_0v_{\sigma}\cdots v_{a\sigma}hv_t$.

The table records, for each type of nonedge, paths of order $r$ avoiding 
the edges that carry $\alpha$.
In each cell a suitable listed path (for a suitable choice of $\sigma$ and $k$) avoids those edges.
\begin{table}[H]
\centering
\[
\begin{array}{|c|c|c|c|}
\hline
& hv_i,\,hv_{i+2d}
& v_jv_{j+1}
&
v_\ell x_\ell,\,v_\ell y_\ell,x_\ell y_\ell
\\
\hline
\quad h,x_0 \quad
&
\begin{array}{c}
hv_{\pm(d-1)}\cdots v_0x_0,\\
hv_{\pm(d-2)}\cdots v_0y_0x_0
\end{array}
&
hv_{\pm(d-1)}\cdots v_0x_0
&
\begin{array}{c}
hv_{d-1}\cdots v_0x_0,\\
hv_{d-2}\cdots v_0y_0x_0
\end{array}
\\ \hline
v_0,v_t
& P^t_\sigma(k,d-1-k)
& P^t_\sigma(k,d-1-k)
& P^t_\sigma(k,d-1-k)
\\ \hline
v_0,x_t
&
\begin{array}{c}
P^t_\sigma(k,d-2-k)\,x_t,\\
P^t_\sigma(k,d-3-k)\,y_tx_t
\end{array}
&
\begin{array}{c}
P^t_\sigma(k,d-2-k)\,x_t,\\
P^t_\sigma(k,d-3-k)\,y_tx_t
\end{array}
&
\begin{array}{c}
P^t_\sigma(k,d-2-k)\,x_t,\\
P^t_\sigma(k,d-3-k)\,y_tx_t,\\
P^{t+\sigma}_\sigma(d-3,0)\,v_tx_t,\\
P^{t+\sigma}_\sigma(d-4,0)\,v_ty_tx_t
\end{array}
\\ \hline
x_0,x_t
&
\begin{array}{c}
x_0P^t_\sigma(k,d-3-k)x_t,\\
x_0y_0P^t_\sigma(k,d-4-k)x_t
\end{array}
&
x_0P^t_\sigma(k,d-3-k)x_t
&
\begin{array}{c}
x_0P^t_\sigma(k,d-3-k)x_t,\\
x_0y_0P^t_\sigma(k,d-4-k)x_t,\\
x_0P^t_\sigma(k,d-4-k)y_tx_t
\end{array}
\\ \hline
\end{array}
\]
\caption{\centering Rows are the endpoints of the nonedge; columns are the edges of color $\alpha$ to be avoided.
Here $\sigma=\pm1$, and $0\le k\le d-s$ where $s\in[4]$ is such that the displayed path contains $P^t_\sigma(k,d-s-k)$.}
\end{table}
All rim and petal colors are unique and disjoint from the spoke colors, so a listed route is rainbow as soon as it uses only one spoke, or two spokes of distinct colors.
For a path with two spokes one may choose an arc of length $L>s$ with $L\not\equiv s\pmod d$, whence $s+1\ge3$ consecutive choices of $k$ give distinct spoke pairs and at least one avoids any prescribed spoke color.
One can check that, for every nonedge $zw$ and every color $\alpha$, at least one of the displayed $z$--$w$ paths of order $r$ is simple, rainbow, and avoids $\alpha$.
This proves the lemma. \qed
\vspace{0.5em}

\noindent {\bf Proof of Theorem~\ref{T15}. }We first show that $(Q_{r,p},\mathcal C)$ is $C_r$-rainbow saturated.
Any cycle that uses a petal edge is a triangle, while the only cycle avoiding $h$ is the rim $C_p$ of length $p\ge4d>r$.
A cycle through $h$ consists of two spokes and a rim path.
It has length $r$ precisely when this rim path has length $d=r-2$, in which case the spokes are $hv_i$ and $hv_{i\pm d}$ and receive the same color.
Thus $(Q_{r,p},\mathcal C)$ contains no rainbow $C_r$.
Lemma~\ref{Lperiodic} and Proposition~\ref{P21} show that it is $C_r$-rainbow saturated.

For an arbitrary $n\geq12d+1$, set
\[
p=d\left\lfloor\frac{n-1}{3d}\right\rfloor,
\qquad s=n-1-3p.
\]
Then $p\geq4d$, $d\mid p$, and $0\leq s<3d\leq p$, so $Q_{r,p}$ contains at least $s$ attached triangles.
Choose $s$ of them and replace each triangle on $\{v_i,x_i,y_i\}$ by a $K_4$ on $v_i$ and three petal vertices; each replacement adds one vertex and three edges.
Denote the resulting graph by $Q_{r,p}^s$.
Lemma~\ref{Lperiodic} can be extended to $Q_{r,p}^s$, which likewise contains no rainbow $C_r$, so Proposition~\ref{P21} yields that $(Q_{r,p}^s,\mathcal C)$ is $C_r$-rainbow saturated.
It has $n$ vertices and
\[
e(Q_{r,p}^s)=5p+3s
=\frac53(n-1)+\frac43s
<\frac53(n-1)+4d.
\]
\qed

\section{Lower bound of $rsat(n,C_r)$}
In this section we consider cycles of length at least five.
Let $n\geq r\geq5$, and let $(G,\mathcal C)$ be a $C_r$-rainbow saturated graph of order $n$.
%Although we have showed a better lower bound for $C_5$, the following several properties hold together for $C_5$ and are useful to help us improve the lower bound.

\begin{prop}\label{P71}
Let $u \in V(G)$ with $d(u)=2$ and $N(u)=\{v,w\}$.
Assume that $d(v) \leq d(w)$. Then exactly one of the following holds.
{\renewcommand{\theenumi}{\arabic{enumi}}
\begin{enumerate}
\item $d(v) \geq 3$;
\item $d(v) = 2$ and $N(v)=\{u,w\}$.
\end{enumerate}}
\end{prop}
\noindent {\bf Proof. }Suppose (1) does not hold, then $d(v)=2$ and $N(v)=\{u,x\}$.
If $x \neq w$, then all $P_r$'s connecting $u,x$ can not avoid $uw$, a contradiction to Corollary \ref{C23}. \qed
\vspace{0.4em}

	For a degree 2 vertex $u$, we say that $u$ is a \textit{good root} (resp. \textit{bad root}) if it satisfies Proposition \ref{P71}(1) (resp. Proposition \ref{P71}(2)).
	Clearly, the number of bad roots in $G$ must be even and any bad root must lie in a triangle together with another bad root.
	We say that $w$ is a \textit{suspension} vertex if there exists bad roots in the neighbourhood of $w$.
	Obviously $w$ is a cut vertex in $G$.
	
\begin{lemma}\label{L72}
If $w$ is a suspension vertex, then the number of bad roots in $N(w)$ is 2.
Moreover, $d(w) \geq 5$ and $N(w)$ has no good roots.
\end{lemma}
\noindent {\bf Proof. }Suppose not, then there exists at least $4$ bad roots in $N(w)$, says $x_1,x_2,y_1,y_2$, with $x_1x_2w$, $y_1y_2w$ being triangles.
Then there exists only $P_3$'s, $P_4$'s, $P_5$'s connecting $x_1,y_1$ and the $P_5$ connecting $x_1,y_1$ can not avoid $x_1x_2$, a contradiction to Corollary \ref{C23}.
Thus, the number of bad roots in $N(w)$ is 2.

	Let $x_1,x_2$ be the bad roots in $N(w)$.
	Suppose $d(w)\leq 4$.
	Let $N(w)\setminus \{x_1,x_2\} = \{y\}$ or $\{y,z\}$.
	Then $x_1,y$ is not connected by any $P_r$, or all $P_r$'s connecting $x_1,y$ can not avoid $wz$ which contradicts Corollary \ref{C23}.
	Hence $d(w) \geq 5$.
	
	Suppose $N(w)$ has a good root $y$, say $N(y)=\{w,z\}$.
	Let $x_1$ be a bad root in $N(w)$.
	Then all $P_r$'s connecting $x_1,y$ can not avoid $zy$, a contradiction to Corollary \ref{C23}.
	Thus, $N(w)$ has no good roots. \qed
\vspace{0.4em}

% 	By Lemma \ref{L72}, we immediately have the following result.
% \vspace{0.2em}

\begin{cor}\label{C73}
If $G$ has no good root, then $e(G) \geq \frac{3}{2}n$.
\end{cor}
\noindent {\bf Proof. }Let $D$ be the set of bad roots in $G$ and $H=G[V(G) \setminus D]$.
By Lemma \ref{L72}, we have $\delta(H) \geq 3$.
Hence $e(G)=e(H)+e(D,V(H))+e(G[D]) \geq \frac{3}{2} |V(H)| + |D| + \frac{1}{2} |D| = \frac{3}{2}n$. \qed
\vspace{0.4em}

\begin{prop}\label{P76}
Let $r\geq6$, and define $F_0=1$, $F_1=3$, and $F_t=3+2F_{t-2}$ for $t\geq2$.
Every $C_r$-rainbow saturated colored graph $(G,\mathcal C)$ of order $n$ satisfies
\[
e(G)\geq
\left(\frac65+\frac{1}{5F_{r-2}}\right)n-O_r(1).
\]
\end{prop}

\noindent {\bf Proof. }Let $A$ be the set of good roots and $D$ the set of bad roots, and put $H=G-D$. 
Clearly $H$ is $C_r$-rainbow saturated, so $diam(H) \leq r-1$.
If $A=\emptyset$, Corollary~\ref{C73} gives $e(G)\geq3n/2$ and we are done.
We may therefore assume that $A\neq\emptyset$.

Set $U=V(H)\setminus A$.
Let $e_0=e(H[U])$ and $f(U)=\sum_{v\in U}\bigl(d_H(v)-3\bigr)$.
Lemma~\ref{L72} gives $d_H(v)\geq3$ for $v\in U$, so $f(U)\geq0$.
The set $A$ is independent and each of its vertices has exactly two neighbors in $U$.
Consequently
\(
2|A|+2e_0=3|U|+f(U), e(H)=2|A|+e_0
\)
and hence
\(
5e(H)-6|H|=2f(U)+e_0.
\)

Take a BFS tree $T$ of $H$ starting at a vertex $x^*\in U$.
First delete every tree edge with both ends in $U$, then at most $e_0$ edges are deleted and a forest remains in which every remaining edge joins $U$ to $A$ since $A$ is independent, so every path alternates between $U$ and $A$.
%A vertex of $U\setminus\{x^*\}$ may now have no parent in this forest, but only because its parent edge in $T$ already lay in $H[U]$ and was deleted in this step; that vertex still has at most $d_H(v)-1$ children in $T$.
Next, at each $v\in U$, among the remaining child-edges of $T$ from $v$ to $A$, retain two of them if $v\neq x^*$ and three if $v=x^*$, or all of them if there are fewer.
Thus at most $d_H(v)-3$ edges are deleted at $v$ in this step, and altogether at most $f(U)$ further edges are deleted.

Every component except the one containing $x^*$ has height at most $r-2$ and alternates between $U$ and $A$.
%We may assume that it starts at a parent in $U$: if the vertex closest to $x^*$ lies in $A$, then that vertex has at most one $U$-child, so the component is no larger than a $U$-started component of height at most $r-2$.
A $U$-vertex has at most two $A$-children, and an $A$-vertex has at most one $U$-child.
Thus such a component has at most $F_{r-2}$ vertices: the extremal sizes satisfy
\[
F_0=1,\qquad F_1=3,\qquad F_t=3+2F_{t-2}.
\]
The component containing $x^*$ has at most $1+3(1+F_{r-3})$ vertices because $x^*$ has at most three $A$-children and below each such child the same alternating recursion applies.
Since deleting $e_0+f(U)$ tree edges creates at most that many additional components,
\[
|H|\leq O_r(1)+F_{r-2}\bigl(e_0+f(U)\bigr).
\]
It follows that
\[
e_0+f(U)\geq\frac{|H|}{F_{r-2}}-O_r(1).
\]
Since $5e(H)-6|H|=2f(U)+e_0$, we obtain
\begin{equation*}
\begin{aligned}
5e(H)
&= 6|H|+2f(U)+e_0 \\
&\geq 6|H|+e_0+f(U) \\
&\geq 6|H|+\frac{|H|}{F_{r-2}}-O_r(1),
\end{aligned}
\end{equation*}
and hence
\[
e(H)\geq
\left(\frac65+\frac{1}{5F_{r-2}}\right)|H|-O_r(1).
\]

Finally, Lemma~\ref{L72} shows that the bad roots occur in disjoint pairs, each pair contributing two vertices and three edges.
Since $3/2\geq6/5+1/(5F_{r-2})$, adding them back preserves the stated inequality. \qed
\vspace{1em}

\noindent{\bf Proof of Theorem \ref{T17}. } By Proposition~\ref{P76} we are done. \qed

% \section{Concluding remarks}

% We determined $\operatorname{rsat}(n,C_4)$ exactly and obtained new bounds for longer cycles.
% The principal remaining problem is to determine the correct linear coefficient for each fixed cycle.

% \begin{prob}\label{C81}
% Does the limit $\lim_{n\to\infty}\operatorname{rsat}(n,C_r)/n$ exist for every fixed $r\geq5$?
% \end{prob}

\section{Concluding remarks}

We determined $\operatorname{rsat}(n,C_4)$ exactly and obtained new bounds for longer cycles. The principal remaining problem is to determine the correct linear coefficient for each fixed cycle. In particular, considerable gaps remain between our upper and lower bounds for $C_5$ and $C_6$, as well as for general cycles of length at least seven.

It would also be interesting to obtain a better structural understanding of sparse $C_r$-rainbow saturated graphs. Our lower-bound argument relies mainly on the behavior of degree-two vertices, while our upper bounds arise from gluing sparse blocks or using periodic constructions. Combining these structural and constructive viewpoints may lead to sharper bounds.

\begin{prob}
	Does the limit
	\(
	\lim_{n\to\infty}\frac{\operatorname{rsat}(n,C_r)}{n}
	\)
	exist for every fixed $r\geq 5$?
\end{prob}

\section*{Acknowledgement}
The authors used GPT-5.6 for grammar and style checking, and to find the construction for $r \geq 6$. 
It was not used to develop, verify, or modify any other proof or mathematical argument in this paper. 
After using this tool, the authors reviewed and edited the content as needed and take full responsibility for the content of the article.

The research of Lu is supported by the National Natural Science Foundation of China (No.~12401445, 12171272). 
The research of He is supported by Beijing Natural Science Foundation (No.~1244047) and China Postdoctoral Science Foundation (No.~2023M740207). 
The research of Xu is supported by China Scholarship Council (No.~202606210184).


\begin{thebibliography}{99}

	\bibitem{BEH}  N. Behague, T. Johnston, S. Letzter, N. Morrison, and S. Ogden, The rainbow saturation number is linear, SIAM J. Discrete Math. {\bf 38} (2024), no.~2, pp. 1239--1249. %rsat
	

	\bibitem{CHA} D. Chakraborti, K. Hendrey, B. Lund and C. Tompkins, Rainbow saturation for complete graphs, SIAM J. Discrete Math. {\bf 38} (2024), no.~1, pp. 1090--1112. %rsat
	
	\bibitem{Chen1} Y. Chen, Minimum $C_5$-saturated graphs, J. Graph Theory, 61(2009), pp. 111–126.
	
	\bibitem{Chen2} Y. Chen, All minimum $C_5$-saturated graphs, J. Graph Theory, 67(2011), pp. 9–26.
	
	\bibitem{survey} B. L. Currie, J. R. Faudree, R. J. Faudree, and J. R. Schmitt, A survey of minimum saturated graphs, Electron. J. Combin., 18 (2011), D519.

	\bibitem{EHM} P. Erd\H os, A. Hajnal and J.W. Moon, A problem in graph theory,  Amer. Math. Monthly  71 (1964), pp. 1107-1110.
	
	
	\bibitem{FI} D. Fisher, K. Fraughnaugh and L. Langley, $P_3$-connected graphs of minimum size, Ars Combin. 47 (1997), pp. 299–306.
	
	\bibitem{FI2} D. Fisher, K. Fraughnaugh and L. Langley, On $C_5$-saturated graphs with minimum size, Graph Theory and Computing (1995), volume 112, pp. 45-48.	
	
 	\bibitem{Fur} Z. F{\"u}redi, Y. Kim. Cycle-saturated graphs with minimum number of edges. J. Graph Theory, 73(2) (2013), pp. 203-215.
	
	\bibitem{Gir} A. Gir\~ao, D.~C. Lewis and K. Popielarz, Rainbow saturation of graphs, J. Graph Theory {\bf 94} (2020), no.~3, pp. 421--444. %rsat
	
	
	\bibitem{Gou}R. Gould, T.  Luczak, J. Schmitt, Constructive upper bounds for cycle-saturated graphs of minimum size, Electron. J. Combin. 13 (2006), D29.
			

	%\bibitem{KT} L. K\'aszonyi and Z. Tuza, Saturated graphs with minimal number of edges, J. Graph Theory, 10 (1986), pp. 203-210.
	
	\bibitem{LAN} Y. Lan, Y. Shi, Y. Wang and J. Zhang, The saturation number of $C_6$, Discrete Math. 348 (2025), 114504.

	\bibitem{MPT} A. Mohammadian, M. Poursoltani and B. Tayfeh-Rezaie, On saturation numbers of complete multipartite graphs and even cycles, arXiv:2506.09767 (2025).

	
	\bibitem{OLL} L. Ollmann. $K_{2,2}$-saturated graphs with a minimal number of edges. Pro. 3rd Southeastern Conference on Combinatorics, Graph and Computing (1972), pp. 367-392.
	
	%\bibitem{MO} G. Moshkovitz and A. Shapira, Exact bounds for some hypergraph saturation problems, J. Combin. Theory Ser. B, 111 (2015), pp. 242-248.
		

    \bibitem{TUZ} Z. Tuza. $C_4$-saturated graphs of minimum size. Acta Univ. Carolin. Math. Phys., 30(2) (1989), pp. 161-167.

%    \bibitem{W1} W. Wessel, \H uber eine Klasse paarer Graphen. I. Beweis einer Vermutung von Erd\H os, Hajnal und Moon, Wiss. Z. Techn. Hochsch. Ilmenau, 12 (1966), pp. 253-256.
%
%	\bibitem{W2} W. Wessel, \H uber eine Klasse paarer Graphen. II. Bestimmung der Minimalgraphen, Wiss. Z. Techn. Hochsch. Ilmenau, 13 (1967), pp. 423-426.
	\bibitem{Zhang} M. Zhang, S. Luo and M. Shigeno, On the Number of Edges in a Minimum $C_6$-Saturated Graph, Graphs
Combin. 31(4)(2015), pp. 1085–1106.



\end{thebibliography}
\end{document}